\documentclass[a4paper,12pt]{article}
\usepackage{latexsym}
\usepackage{amsfonts}
\newtheorem{theorem}{Theorem}[section]
\newtheorem{corollary}[theorem]{Corollary}
\newtheorem{definition}{Definition}
\newtheorem{example}[theorem]{Example}
\newtheorem{remark}[theorem]{Remark}
\newtheorem{lemma}[theorem]{Lemma}
\newtheorem{proposition}[theorem]{Proposition}
\title{Idealization of Ganster--Reilly decomposition
theorems\thanks{1991 Math.\ Subject Classification --- Primary:
54C08, 54A20; Secondary: 54A05, 54C10. \protect\newline Key words
and phrases --- decomposition of continuity, topological ideal,
$\cal I$-open set, preopen, locally closed, pre-$\cal I$-open
set, $\cal I$-continuous function. \protect\newline Research
supported partially by the Ella and Georg Ehrnrooth Foundation
at Merita Bank, Finland.}}
\author{Julian Dontchev\\Department of Mathematics\\University
of Helsinki\\PL 4, Yliopistonkatu 5\\00014 Helsinki\\Finland}
\date{}
\begin{document}
\baselineskip=20pt plus 1pt minus 1pt
\newcommand{\pio}{pre-$\cal I$-open}
\newcommand{\pic}{pre-$\cal I$-continuous}
\newcommand{\sic}{$\star$-$\cal I$-continuous}
\newcommand{\xti}{$(X,\tau,{\cal I})$}
\newcommand{\fxy}{$f \colon (X,\tau) \rightarrow (Y,\sigma)$}
\maketitle
\begin{abstract} 
In 1990, Ganster and Reilly \cite{GR1} proved that a function
\fxy\ is continuous if and only if it is precontinuous and
$LC$-continuous. In this paper we extend their decomposition of
continuity in terms of ideals. We show that a function
$f \colon (X,\tau,{\cal I}) \rightarrow (Y,\sigma)$ is
continuous if and only if it is $\cal I$-continuous and
$\cal I$-$LC$-continuous. We also provide a decomposition of
$\cal I$-continuity.
\end{abstract}

\section{Introduction to topological ideals}\label{s1}

In \cite{GR1,GR2,GR3}, Ganster and Reilly gave several new
decompositions of continuity.

Let $A$ be a subset of a topological space $(X,\tau)$.
Following Kronheimer \cite{Kr1}, we call the interior of the
closure of $A$, denoted by $A^+$, the {\em consolidation} of $A$.
Sets included in their consolidation are called {\em preopen}
or {\em locally dense} \cite{CM1}. If $A$ is the intersection of
an open and a closed (resp.\ regular closed) set, then $A$ is
called {\em locally closed} (resp.\ {\em $\cal A$-set}
\cite{Tong1}). A function \fxy\ is called {\em precontinuous}
(resp.\ {\em $LC$-continuous} \cite{GR4}, {\em $\cal
A$-continuous} \cite{Tong1}) if the preimage of every open set
is preopen (resp.\ locally closed, $\cal A$-set). The following
theorem is due to Ganster and Reilly \cite[Theorem 4 (iv) and
(v)]{GR1}. 

\begin{theorem}\label{grt1}
{\em \cite{GR1}} For a function \fxy\ the following conditions
are equivalent:

{\rm (1)} $f$ is continuous.

{\rm (2)} $f$ is precontinuous and $\cal A$-continuous.

{\rm (3)} $f$ is precontinuous and $LC$-continuous.
\end{theorem}

The aim of this paper is to present an idealized version of the
Ganster--Reilly decomposition theorem.

A nonempty collection $\cal I$ of subsets on a topological space
$(X,\tau)$ is called an {\em ideal} on $X$ if it satisfies the
following two conditions:

(1) If $A \in {\cal I}$ and $B \subseteq A$, then $B \in {\cal
I}$ (heredity).

(2) If $A \in {\cal I}$ and $B \in {\cal I}$, then $A \cup B \in
{\cal I}$ (finite additivity).

A {\em $\sigma$-ideal} on a topological space $(X,\tau)$ is an
ideal which satisfies:

(3) If $\{ A_{i} \colon i = 1,2,3, \ldots \} \subseteq {\cal I}$,
then $\bigcup \{ A_{i} \colon i = 1,2,3, \ldots \} \in {\cal I}$
(countable additivity).

If $X \not\in {\cal I}$, then $\cal I$ is called a {\em proper}
ideal. The collection of the complements of all elements of a
proper ideal is a filter, hence proper ideals are sometimes
called {\em dual filters}.

The following collections form important ideals on a topological
space $(X,\tau)$: the ideal of all finite sets $\cal F$, the
$\sigma$-ideal of all countable sets $\cal C$, the ideal of all
closed and discrete sets $\cal CD$, the ideal of all nowhere
dense sets $\cal N$, the $\sigma$-ideal of all meager sets $\cal
M$, the ideal of all scattered sets $\cal S$ (here $X$ must be
$T_0$) and the $\sigma$-ideal of all Lebesgue null sets $\cal L$
(here $X$ is the real line).

An {\em ideal topological space} is a topological space
$(X,\tau)$ with an ideal $\cal I$ on $X$ and is denoted by
$(X,\tau,{\cal I})$. For a subset $A \subseteq X$, $A^{*}({\cal
I}) = \{x \in X \colon$ for every $U \in {\tau} (x), U \cap A
\not\in {\cal I} \}$ is called the {\em local function\/} of $A$
with respect to $\cal I$ and $\tau$ \cite{JH1,K1}. We simply
write $A^{*}$ instead of $A^{*}({\cal I})$ in case there is no
chance for confusion. Note that often $X^*$ is a proper subset
of $X$. The hypothesis $X = X^*$ was used by Hayashi in
\cite{H1}, while the hypothesis $\tau \cap {\cal I} = \emptyset$
was used by Samuels in \cite{S1}. In fact, those two conditions
are equivalent \cite[Theorem 6.1]{JH1} and we call the ideal
topological spaces which satisfy this hypothesis {\em
Hayashi-Samuels spaces}. Note that $(X,\tau,\{\emptyset\})$
and $(X,\tau,{\cal N})$ are always Hayashi-Samuels spaces; also
$({\mathbb R},\tau,{\cal F})$ is a Hayashi-Samuels space, where
$\tau$ denotes the usual topology on the real line $\mathbb R$.

For every ideal topological space $(X,\tau,{\cal I})$, there
exists a topology $\tau^{*}({\cal I})$, finer than $\tau$,
generated by the base $\beta({\cal I},\tau) = \{ U \setminus I
\colon U \in \tau$ and $I \in {\cal I}\}$. In general,
${\beta}({\cal I},\tau)$ is not always a topology \cite{JH1}.
When there is no chance for confusion, $\tau^{*}({\cal I})$ is
denoted by $\tau^{*}$. Observe additionally that ${\rm Cl}^{*}(A)
= A \cup A^{*}$ defines a Kuratowski closure operator for (the
same topology) $\tau^{*}({\cal I})$.

Recall that $A \subseteq (X,\tau,{\cal I})$ is called {\em
$\star$-dense-in-itself} \cite{H1} (resp.\ {\em
$\tau^{*}$-closed} \cite{JH1}, {\em $\star$-perfect}
\cite{H1}) if $A \subseteq A^{*}$ (resp.\ $A^{*} \subseteq A$,
$A = A^{*}$).

It is interesting to note that $A^{*}({\cal I})$ is a
generalization of closure points, $\omega$-accumulation points
and condensation points. Recall that the set of all
$\omega$-accumulation points of subset $A$ of a topological space
$(X,\tau)$ is $A^{\omega} = \{ x \in X \colon U \cap A$ is
infinite for every $U \in {\cal N}(x) \}$. The set of all
condensation points of $A$ is $A^{k} = \{ x \in X \colon U
\cap A$ is uncountable for every $U \in {\cal N}(x) \}$. It is
easily seen that ${\rm Cl} (A) = A^{*}(\{ \emptyset \})$,
$A^{\omega} = A^{*}({\cal F})$ and $A^{k} = A^{*}({\cal C})$.
Note here that in $T_1$-spaces the concepts of
$\omega$-accumulation points and limit points coincide.

In 1990, D. Jankovi\'{c} and T.R. Hamlett introduced the notion
of $\cal I$-open sets in ideal topological spaces. Given an ideal
topological
space $(X,\tau,{\cal I})$ and $A \subseteq X$, $A$ is said to be
{\em $\cal I$-open} \cite{JH2} if $A \subseteq {\rm Int}(A^{*})$.
We denote by $IO(X,\tau,{\cal I}) = \{A \subseteq X: A \subseteq
{\rm Int}(A^{*})\}$ or simply write $IO(X,\tau)$ or $IO(X)$ when
there is no chance for confusion with the ideal. A subset $F
\subseteq (X,\tau,{\cal I})$ is called \cite{A1} {\em $\cal
I$-closed} if its complement is $\cal I$-open. Note that $X$ need
not be an ${\cal I}$-open subset. Thus, not only are ${\cal
I}$-open and $\tau^{*}$-open sets are different concepts, but the
former do not give a topology. In the extreme case when ${\cal
I}$ is the maximal ideal of all subsets of $X$, only the void
subset is ${\cal I}$-open.

A function $f:(X,\tau,{\cal J}_{1}) \rightarrow
(Y,\sigma,{\cal J}_{2})$ is said to be {\em $\cal I$-continuous}
(resp.\ {\em $\cal I$-open}, {\em $\cal I$-closed}) if for every
$V \in \sigma$ (resp.\ $U \in \tau$, $U$ closed in $X$),
$f^{-1}(V) \in IO(X,\tau)$ (resp.\ $f(U) \in IO(X,\tau)$, $f(U)$
is $\cal I$-closed). The definitions are due to Monsef {\em et
al.\/} \cite{A1}.

In \cite{Nj1}, a topology ${\tau}^{\alpha}$ has
been introduced by defining its open sets to be the
$\alpha$-sets, that is the sets $A \subseteq X$ with $A \subset
{\rm Int}({\rm Cl}({\rm Int}(A)))$. Observe that $\tau^{\alpha}
= \tau^{*}({\cal N})$.

\section{Pre-$\cal I$-open sets}\label{s2}

\begin{definition}\label{d1}
{\em A subset of an ideal topological space $(X,\tau,{\cal I})$
is called {\em \pio\ } if $A \subseteq {\rm Int}({\rm
Cl}^{*}(A))$.}
\end{definition}

We denote by $PIO(X,\tau,{\cal I})$ the family of all \pio\
subsets of $(X,\tau,{\cal I})$ or simply write $PIO(X,\tau)$ or
$PIO(X)$ when there is no chance for confusion with the ideal.
We call a subset $A \subseteq (X,\tau,{\cal I})$ {\em pre-$\cal
I$-closed} if its complement is \pio.

Although $\cal I$-openness and openness are independent concepts
\cite[Examples 2.1 and 2.2]{A1}, pre-$\cal I$-openness is related
to both of them as the following two results show.

\begin{proposition}\label{t1}
Every $\cal I$-open set is \pio.
\end{proposition}

{\em Proof.} Let $(X,\tau,{\cal I})$ be an ideal topological space and
let $A \subseteq X$ be $\cal I$-open. Then $A \subseteq {\rm Int}
(A^{*}) \subseteq {\rm Int} (A^{*} \cup A) = {\rm Int}({\rm
Cl}^{*} (A))$. $\Box$

\begin{proposition}\label{t2}
Every open set is \pio.
\end{proposition}

{\em Proof.} Let $A \subseteq (X,\tau,{\cal I})$ be open. Then
$A \subseteq {\rm Int} A \subseteq {\rm Int} (A^{*} \cup A) =
{\rm Int}({\rm Cl}^{*} (A))$. $\Box$

The converse in the proposition above is not necessarily true as
shown by the following two examples.

\begin{example}\label{e1}
{\em A \pio\ set, even an open set, need not be $\cal I$-open.
Let $X = \{ a,b,c,d \}$, $\tau = \{ \emptyset, \{ a,c \}, \{ d
\}, \{ a,c,d \}, X \}$ and ${\cal I} = \{ \emptyset, \{ c \}, \{
d \}, \{ c,d \} \}$. Set $A = \{ a,c,d \}$. Then $A \in \tau$ and
hence $A \in PIO(X)$ but $A \not\in IO(X)$ \cite[Example
2.2]{A1}.}
\end{example}

\begin{example}\label{e2}
{\em Let $(X,\tau)$ be the real line with the usual topology and
let ${\cal F}$ be as mentioned before the ideal of all finite
subsets of $X$. Let ${\mathbb Q}$ be the set of all rationals.
Since ${\mathbb Q}^{*}({\cal F}) = X$, then ${\mathbb Q}$ is
\pio\ (even $\cal I$-open). But clearly ${\mathbb Q} \not\in
\tau$.}
\end{example}

Our next two results together with Propsoition~\ref{t1} and
Propsoition~\ref{t2} shows that the class of \pio\ sets is
properly placed between the classes of $\cal I$-open and preopen
sets as well as between the classes of open and preopen sets.

\begin{proposition}\label{t3}
Every \pio\ set is preopen.
\end{proposition}

{\em Proof.} Let $(X,\tau,{\cal I})$ be an ideal topological
space and let $A \in PIO(X)$. Then $A \subseteq {\rm Int}({\rm Cl}^{*}
(A)) = {\rm Int} (A^{*} \cup A) \subseteq {\rm Int} ({\rm Cl} (A)
\cup A) = {\rm Int} ({\rm Cl} (A))$. $\Box$

\begin{example}\label{e3}
{\em A preopen set need not be \pio. Every singleton (for
example) in an indiscrete topological space with cardinality at
least two is preopen but if we set $\cal I$ to be the maximal
ideal, i.e., ${\cal I} = {\cal P}(X)$, then it is easy to see that
none of the singletons is \pio.}
\end{example}

\begin{proposition}\label{t4}
For an ideal topological space $(X,\tau,{\cal I})$ and $A \subseteq X$
we have:

{\rm (i)} If ${\cal I} = \emptyset$, then $A$ is \pio\ if and
only if $A$ is preopen.

{\rm (ii)} If ${\cal I} = {\cal P}(X)$, then $A$ is \pio\ if and
only if $A \in \tau$.

{\rm (iii)} If ${\cal I} = {\cal N}$, then $A$ is \pio\ if
and only if $A$ is preopen.
\end{proposition}

{\em Proof.} (i) Necessity is shown in Proposition~\ref{t3}. For
sufficiency note that in the case of the minimal ideal $A^{*} =
{\rm Cl} (A)$.

(ii) Necessity: If $A \in PIO(X)$, then $A \subseteq {\rm
Int}({\rm Cl}^{*} (A)) = {\rm Int}(A \cup A^{*}) = {\rm Int}(A
\cup \emptyset) = {\rm Int} A$. Sufficiency is given in
Proposition~\ref{t2}.

(iii) By Proposition~\ref{t3} we need to show only sufficiency.
Note that the local function of $A$ with respect to ${\cal
N}$ and $\tau$ can be given explicitly \cite{V2}. We have:
$$A^{*}({\cal N}) = {\rm Cl}({\rm Int}({\rm Cl}(A))).$$ Thus
$A$ is \pio\ if and only if $A \subseteq {\rm Int}(A \cup {\rm
Cl}({\rm Int}({\rm Cl}(A))))$. Assume that $A$ is preopen. Since
always ${\rm Int}({\rm Cl}(A)) \subseteq A \cup {\rm Cl}({\rm
Int}({\rm Cl}(A)))$, then $A \subseteq {\rm Int}(A \cup {\rm
Cl}({\rm Int}({\rm Cl}(A)))) = {\rm Int}(A \cup A^{*}({\cal N}))
= {\rm Int}({\rm Cl}^{*}(A))$ or equivalently $A$ is \pio. $\Box$

The intersection of even two \pio\ sets need not be \pio\ as
shown in the following example.

\begin{example}\label{e4}
{\em Let $X = \{ a,b,c \}$, $\tau = \{ \emptyset, \{ a,b \}, X
\}$ and ${\cal I} = \{ \emptyset, \{ c \} \}$. Set $A = \{ a,c
\}$ and $B = \{ b,c \}$. Since $A^{*} = B^{*} = X$, then both $A$
and $B$ are \pio. But on the other hand $A \cap B = \{ c \}
\not\in PIO(X)$.}
\end{example}

\begin{lemma}\label{l1}
{\em \cite[Theorem 2.3 (g)]{JH1}} Let $(X,\tau,{\cal I})$ be an
ideal topological space and let $A \subseteq X$. Then $U \in \tau
\Rightarrow U \cap A^{*} = U \cap (U \cap A)^{*} \subseteq (U
\cap A)^{*}$. $\Box$
\end{lemma}

\begin{proposition}\label{t5}
Let $(X,\tau,{\cal I})$ be an ideal topological space with
$\triangle$ an arbitrary index set. Then:

{\rm (i)} If $\{ A_{\alpha} \colon \alpha \in \triangle \}
\subseteq PIO(X)$, then $\cup \{ A_{\alpha} \colon \alpha \in
\triangle \} \in PIO(X)$.

{\rm (ii)} If $A \in PIO(X)$ and $U \in \tau$, then $A \cap U \in
PIO(X)$.

{\rm (iii)} If $A \in PIO(X)$ and $B \in \tau^{\alpha}$, then $A
\cap B \in PO(X)$.

{\rm (iv)} If $A \in PIO(X)$ and $B \in SO(X)$, then $A \cap B
\in SO(A)$.

{\rm (v)} If $A \in PIO(X)$ and $B \in SO(X)$, then $A \cap B
\in PO(B)$.

\end{proposition}

{\em Proof.} (i) Since $\{ A_{\alpha} \colon \alpha \in \triangle
\} \subseteq PIO(X)$, then $A_{\alpha} \subseteq {\rm Int}({\rm
Cl}^{*}(A_{\alpha}))$ for every $\alpha \in \triangle$. Thus
$\cup_{\alpha \in \triangle} A_{\alpha} \subseteq \cup_{\alpha
\in \triangle} {\rm Int}({\rm Cl}^{*}(A_{\alpha})) \subseteq {\rm
Int} (\cup_{\alpha \in \triangle} {\rm Cl}^{*}(A_{\alpha})) =
{\rm Int} (\cup_{\alpha \in \triangle} (A_{\alpha}^{*} \cup
A_{\alpha})) = {\rm Int} ((\cup_{\alpha \in \triangle}
A_{\alpha}^{*}) \cup (\cup_{\alpha \in \triangle} A_{\alpha}))
\subseteq {\rm Int} ((\cup_{\alpha \in \triangle} A_{\alpha})^{*}
\cup (\cup_{\alpha \in \triangle} A_{\alpha})) = {\rm Int}({\rm
Cl}^{*} (\cup_{\alpha \in \triangle} A_{\alpha}))$.

(ii) By assumption $A \subseteq {\rm Int}({\rm Cl}^{*}(A))$ and
$U \subseteq {\rm Int} (U)$. Thus applying Lemma~\ref{l1}, $A \cap
U \subseteq {\rm Int}({\rm Cl}^{*}(A)) \cap {\rm Int} (U) \subset
{\rm Int} ({\rm Cl}^{*}(A) \cap U) = {\rm Int} ((A^{*} \cup A)
\cap U) = {\rm Int} ((A^{*} \cap U) \cup (A \cap U)) \subseteq
{\rm Int} ((A \cap U)^{*} \cup (A \cap U)) = {\rm Int}({\rm
Cl}^{*} (A \cap U))$.

(iii) Since the intersection of a preopen set and an $\alpha$-set
is always a preopen set, then the claim is clear due to
Propsoition~\ref{t3}.

(iv) and (v) It is proved in \cite{Noiri1} that the intersection
of a preopen and a semi-open set is a preopen subset of the
semi-open set and a semi-open subset of the preopen set. Thus the
claim follows from Proposition~\ref{t3}. $\Box$

\begin{corollary}\label{c1}
{\rm (i)} The intersection of an arbitrary family of pre-$\cal
I$-closed sets is a pre-$\cal I$-closed set.

{\rm (ii)} The union of a pre-$\cal I$-closed set and a closed
set is pre-$\cal I$-closed. $\Box$
\end{corollary}

Recall that $(X,\tau)$ is called {\em submaximal} if every dense
subset of $X$ is open.

\begin{lemma}
{\em \cite[Lemma 5]{MR1}} If $(X,\tau)$ is submaximal, then
$PO(X,\tau) = \tau$. $\Box$
\end{lemma}

\begin{corollary}
If $(X,\tau)$ is submaximal, then for any ideal $\cal I$ on $X$,
$\tau = PIO(X)$. $\Box$
\end{corollary}

\begin{remark}
{\em By Proposition~\ref{t5}, the intersection of a \pio\ set and
an open set is \pio. However, the intersection of a \pio\ set and
an $\cal I$-open set is not necessarily \pio, since in
Example~\ref{e4} $\{ c \} = A \cap B$ is not \pio, although $A$
is \pio\ (even $\cal I$-open) and $B$ is $\cal I$-open.}
\end{remark}

\begin{remark}
{\em (i) In an ideal topological space $(X,\tau,{\cal I})$, the subset
$X$ need not always be $\cal I$-open. However, $X$ is always \pio.

(ii) If $A \subseteq (X,\tau,{\cal I})$ is $\star$-perfect, then
$A \in \tau$ if and only if $A \in IO(X)$ if and only if $A \in
PIO(X)$.}
\end{remark}

{\bf Problem.} The class of ideal topological spaces \xti\ with
$PIO(X,\tau,{\cal I}) \subseteq \tau^{*}({\cal I})$ is probably
of some interest. Call these spaces {\em $\cal I$-strongly
irresolvable}. It is not difficult to observe that in the trivial
case ${\cal I} = \{ \emptyset \}$, we have the class of strongly
irresolvable spaces which were introduced in 1991 by Foran and
Liebnitz \cite{FL1}. Note also that in the case of the maximal
ideal ${\cal P}(X)$, every ideal topological space is ${\cal
P}(X)$-strongly irresolvable. It is the author's belief that
further study of $\cal I$-strongly irresolvable spaces is
worthwhile.

\section{A decomposition of $\cal I$-continuity}\label{s3}

\begin{definition}\label{d2}
{\em A function $f \colon (X,\tau,{\cal I}) \rightarrow
(Y,\sigma)$ is called {\em \pic\ } if for every $V \in \sigma$,
$f^{-1} (V) \in PIO(X,\tau)$.}
\end{definition}

In the notion of Proposition~\ref{t2} we have the following
result:

\begin{proposition}\label{tt1}
Every continuous function $f \colon (X,\tau,{\cal I}) \rightarrow
(Y,\sigma)$ is \pic. $\Box$ 
\end{proposition}

The converse is not true in general as shown in the following
example.

\begin{example}\label{ee1}
{\em Consider first the classical Dirichlet function $f \colon
{\mathbb R} \rightarrow {\mathbb R}$:

\[ f(x) = \left\{ \begin{array}{ll} 1, &
\mbox{$x \in {\mathbb Q},$} \\ 0, & \mbox{otherwise.} \end{array}
\right. \]

Let ${\cal F}$ be the ideal of all finite subsets of ${\mathbb
R}$. The Dirichlet function $f \colon ({\mathbb R},\tau,{\cal F})
\rightarrow ({\mathbb R},\tau)$ is \pic, since every point of
${\mathbb R}$ belongs to the local function of the rationals with
respect to ${\cal F}$ and $\tau$ as well as to the local function
of the irrationals. Hence $f$ is even $\cal I$-continuous. But
on the other hand the Dirichlet function is not continuous at any
point of its domain.}
\end{example}

Due to Proposition~\ref{t1} we have the next result:

\begin{proposition}\label{tt2}
Every $\cal I$-continuous function $f \colon (X,\tau,{\cal I})
\rightarrow (Y,\sigma)$ is \pic. $\Box$ 
\end{proposition}

The reverse is again not true as the following example shows.

\begin{example}\label{ee2}
{\em Let $(X,\tau,{\cal I})$ be the space from Example~\ref{e1}
and let $\sigma = \{ \emptyset, \{ a,c,d \}, X \}$. Then the
identity function $f \colon (X,\tau,{\cal I}) \rightarrow
(X,\sigma)$ is \pic\ but not $\cal I$-continuous.}
\end{example}

From Proposition~\ref{t3} we have:

\begin{proposition}\label{tt3}
Every \pic\ function $f \colon (X,\tau,{\cal I}) \rightarrow
(Y,\sigma)$ is precontinuous. $\Box$ 
\end{proposition}

\begin{example}\label{ee3}
{\em A precontinuous function need not be \pic. Let $(X,\tau)$
be the real line with the indiscrete topology and $(Y,\sigma)$
the real line with the usual topology. The identity function
$f \colon (X,\tau,{\cal P}(X)) \rightarrow (Y,\sigma)$ is
precontinuous but not \pic.}
\end{example}

\begin{proposition}\label{tt4}
For a function $f \colon (X,\tau,{\cal I}) \rightarrow
(Y,\sigma)$ the following conditions are equivalent:

{\rm (1)} $f$ is \pic.

{\rm (2)} For each $x \in X$ and each $V \in \sigma$ containing
$f(x)$, there exists $W \in PIO(X)$ containing $x$ such that
$f(W) \subseteq V$.

{\rm (3)} For each $x \in X$ and each $V \in \sigma$ containing
$f(x)$, ${\rm Cl}^{*}(f^{-1}(V))$ is a neighborhood of $x$.

{\rm (4)} The inverse image of each closed set in $(Y,\sigma)$
is pre-$\cal I$-closed.
\end{proposition}

{\em Proof.} (1) $\Rightarrow$ (2) Let $x \in X$ and let $V \in
\sigma$ such that $f(x) \in V$. Set $W = f^{-1} (V)$. By (1), $W$
is \pio\ and clearly $x \in W$ and $f(W) \subseteq V$.

(2) $\Rightarrow$ (3) Since $V \in \sigma$ and $f(x) \in V$, then
by (2) there exists $W \in PIO(X)$ containing $x$ such that $f(W)
\subseteq V$. Thus, $x \in W \subseteq {\rm Int}({\rm Cl}^{*}(W))
\subseteq {\rm Int}({\rm Cl}^{*}(f^{-1}(V))) \subseteq {\rm
Cl}^{*}(f^{-1}(V))$. Hence, ${\rm Cl}^{*}(f^{-1}(V))$ is a
neighborhood of $x$.

(3) $\Rightarrow$ (1) and (1) $\Leftrightarrow$ (4) are obvious.
$\Box$

The composition of two \pic\ functions need not be always \pic\
as the following example shows.

\begin{example}\label{ee4}
{\em Let ${\mathbb R}$ be again the real line and $\tau$ the
usual topology. Note that the identity function $g \colon
({\mathbb R},\tau,{\cal P}(X)) \rightarrow ({\mathbb
R},\tau,{\cal F})$ is \pic\ and also the Dirichlet function $f
\colon ({\mathbb R},\tau,{\cal F}) \rightarrow ({\mathbb
R},\tau)$ is \pic\ (Example~\ref{ee1}). But their composition $(f
\circ g) \colon ({\mathbb R},\tau,{\cal P}(X)) \rightarrow
({\mathbb R},\sigma)$ is not \pic, since (for example) $f^{-1}
\{(0,2)\} = {\mathbb Q} \not\in PIO({\mathbb R},\tau,{\cal
P}(X))$.}
\end{example}

However the following result holds.

\begin{proposition}\label{tt5}
Let $f \colon (X,\tau,{\cal I}) \rightarrow (Y,\sigma,{\cal J})$
and $g \colon (Y,\sigma,{\cal J}) \rightarrow (Z,\upsilon)$ be
two functions, where $\cal I$ and $\cal J$ are ideals on
$X$ and $Y$ respectively. Then:

{\rm (i)} $g \circ f$ is \pic, if $f$ is \pic\ and $g$ is
continuous.

{\rm (ii)} $g \circ f$ is precontinuous, if $g$ is continuous and
$f$ is \pic.
\end{proposition}

{\em Proof.} Obvious. $\Box$

Hayashi \cite{H1} defined a set $A$ to be $\star$-dense-in-itself
if $A \subseteq A^{*}({\cal I})$. We say that a function $f
\colon (X,\tau,{\cal I}) \rightarrow (Y,\sigma)$ is {\em \sic\
} if the preimage of every open set in $(Y,\sigma)$ is
$\star$-dense-in-itself in $(X,\tau,{\cal I})$. In what follows,
we try do decompose $\cal I$-continuity but before that we will
give a decomposition of $\cal I$-openness. Our next two examples
(the ones after Proposition~\ref{tt6} and Proposition~\ref{tt7})
will show that pre-$\cal I$-continuity and $\star$-$\cal
I$-continuity are independent concepts.

\begin{proposition}\label{tt6}
For a subset $A \subseteq (X,\tau,{\cal I})$ the following
conditions are equivalent:

{\rm (1)} $A$ is $\cal I$-open.

{\rm (2)} $A$ is \pio\ and $\star$-dense-in-itself.
\end{proposition}

{\em Proof.} (1) By Proposition~\ref{t1}, every $\cal I$-open set
is \pio. On the other hand $A \subseteq {\rm Int}(A^{*}) \subset
A^{*}$, which shows that $A$ is $\star$-dense-in-itself.

(2) $\Rightarrow$ (1) By assumption $A \subseteq {\rm Int}({\rm
Cl}^{*}(A)) = {\rm Int}(A^{*} \cup A) = {\rm Int}(A^{*})$ or
equivalently $A$ is $\cal I$-open. $\Box$

Thus we have the following decomposition of $\cal I$-continuity:

\begin{theorem}\label{tt7}
For a function $f \colon (X,\tau{\cal I}) \rightarrow (Y,\sigma)$
the following conditions are equivalent:

{\rm (1)} $f$ is $\cal I$-continuous.

{\rm (2)} $f$ is \pic\ and \sic. $\Box$
\end{theorem}

\begin{example}\label{ee4a}
{\em The identity function $f \colon ({\mathbb R},\tau,{\cal
P}(X)) \rightarrow ({\mathbb R},\tau)$, where $\tau$ stands for
the usual topology on the real line is \pic\ as mentioned in
Example~\ref{ee4} but not \sic, since the local function of every
subset of ${\mathbb R}$ with respect to ${\cal P}(X)$ and $\tau$
coincides with the void set.}
\end{example}

\begin{example}\label{ee4b}
{\em Note that in the case of the minimal ideal every function
is \sic, since the local function of every set coincides with its
closure. But since not every function is precontinuous, then
$\star$-$\cal I$-continuity does not always imply pre-$\cal
I$-continuity.}
\end{example}

\begin{remark}
{\em Of course a very appropriate example would be the
construction of a space with a fixed ideal on it and finding
topologies on the space such that certain functions would show
the independence of pre-$\cal I$-continuity and $\star$-$\cal
I$-continuity as well as the fact that they are both weaker
than $\cal I$-continuity. Such an example is the following: Let
$X = \{ a,b,c \}$, ${\cal I} = \{ \emptyset, \{ c \} \}$, $\tau
= \{ \emptyset, \{ b \}, X \}$, $\sigma = \{ \emptyset, \{ c \},
X \}$, $\nu = \{ \emptyset, \{ a \}, X \}$. The identity function
$f \colon (X,\tau,{\cal I}) \rightarrow (X,\nu,{\cal I})$ is
\sic\ but neither $\cal I$-continuous nor \pic. On the other hand
the identity function $g \colon (X,\sigma,{\cal I}) \rightarrow
(X,\sigma,{\cal I})$ is \pic\ but neither
$\cal I$-continuous nor \sic.}
\end{remark}

In the case when ${\cal N}$ is the ideal of all nowhere dense
subsets precontinuity coincides with pre-$\cal I$-continuity,
while $\beta$-continuity is equivalent to
$\star$-$\cal I$-continuity due to Proposition~\ref{t4}. Recall
that a function $f \colon (X,\tau) \rightarrow (Y,\sigma)$ is
called {\em $\beta$-continuous} (or sometimes {\em
semi-precontinuous}) if the preimage of every open set in
$(Y,\sigma)$ is $\beta$-open in $(X,\tau)$, where a set $A$
is called $\beta$-open if $A \subseteq {\rm Cl}({\rm
Int}({\rm Cl}(A)))$. It is clear, since every preopen set is
$\beta$-open but not vice versa, that the family of all \pio\
subsets of an ideal topological space $(X,\tau,{\cal I})$ is a proper
subset of the family of all $\beta$-open sets.

Consider next the ideal of all meager subsets. Recall that a set
is {\em meager} if it is a countable union of nowhere dense
sets. Meager sets are called often sets of {\em first category}.
If a set is
not meager it is said to be of {\em second category}. The points
of second category of $A$ are the points of $A^{*}({\cal M})$
\cite{K1}. In 1922 Blumberg \cite{Blu1} called a point $x$ of a
space $(X,\tau)$ {\em inexhaustibly approached} by $A \subseteq
X$ if $x \in A^{*}({\cal M})$. If we call the set $A$ {\em
inexhaustibly approached} when every point of $A$ is
inexhaustibly approached by $A$, then clearly a function is
$\star$-${\cal M}$-continuous if and only if the inverse
image of every open set is inexhaustibly approached.

\section{Idealized Ganster--Reilly decomposition
theorem}\label{s4}

A subset $A$ of an ideal topological space \xti\ is called {\em
$\cal I$-locally closed} if $A = U \cap V$, where $U \in \tau$
is $V$ is $\star$-perfect. Note that in the case of the minimal
ideal, $\cal I$-locally closed is equivalent to locally closed,
while $\cal N$-locally closed is equivalent to the Tong's notion
of an $\cal A$-set from \cite{Tong1}.

\begin{proposition}\label{tt42}
For a subset $A \subseteq (X,\tau,{\cal I})$ of a Hayashi-Samuels
space the following conditions are equivalent:

{\rm (1)} $A$ is open.

{\rm (2)} $A$ is \pio\ and $\cal I$-$LC$-continuous.
\end{proposition}

{\em Proof.} (1) $\Rightarrow$ (2) The first part is
Proposition~\ref{t2}. For the second part, note that $A = A \cap
X$, where $A \in \tau$ and $X$ is $\star$-perfect.

(2) $\Rightarrow$ (1) By assumption $A \subseteq {\rm Int}({\rm
Cl}^{*}(A)) = {\rm Int}({\rm Cl}^{*}(U \cap V))$, where $U \in
\tau$ and $V$ is $\star$-perfect. Hence, $A = U \cap A \subseteq
U \cap ({\rm Int}({\rm Cl}^{*}(U)) \cap {\rm Int}({\rm
Cl}^{*}(V))) = U \cap {\rm Int}(V \cup V^{*}) = {\rm Int}(U) \cap
{\rm Int}(V) = {\rm Int} (U \cap V) = {\rm Int}(A)$. This is
shows that $A \in \tau$. $\Box$

\begin{definition}\label{dd41}
{\em A function $f \colon (X,\tau,{\cal I}) \rightarrow
(Y,\sigma)$ is called {\em $\cal I$-$LC$-continuous} if for every
$V \in \sigma$, $f^{-1} (V)$ is $\cal I$-$LC$-closed.}
\end{definition}

\begin{proposition}\label{tt41}
Let \xti\ be a Hayashi-Samuels space. Then, every continuous
function
$f \colon (X,\tau,{\cal I}) \rightarrow (Y,\sigma)$ is $\cal
I$-$LC$-continuous. $\Box$ 
\end{proposition}

The converse is not true in general, since in the case of the
minial ideal \xti\ is a Hayashi-Samuels space but (usual)
$LC$-continuous functions need not be $LC$-continuous \cite{GR4}.

Now, in the notion of Proposition~\ref{tt42}, we have the
following idealized decomposition of continuity:

\begin{theorem}\label{tt43}
Let \xti\ be a Hayashi-Samuels space. For a function $f \colon
(X,\tau,{\cal I}) \rightarrow (Y,\sigma)$ the following
conditions
are equivalent:

{\rm (1)} $f$ is continuous.

{\rm (2)} $f$ is \pic\ and $\cal I$-$LC$-continuous. $\Box$
\end{theorem}

\begin{remark}
{\em From the particular cases $\cal I = \{ \emptyset \}$ and
${\cal I} = {\cal N}$ in Theorem~\ref{tt43} we derive the
well-known Ganster--Reilly decomposition Theorem~\ref{grt1}.}
\end{remark}

\
E-mail: {\tt dontchev@cc.helsinki.fi}, {\tt
dontchev@e-math.ams.org}
\ 
\end{document}